\documentstyle[12pt,fleqn,leqno]{article}
\title{{\LARGE\bf On a theorem of Shapiro}\\}
\author{%
Saka\'e Fuchino, Saharon Shelah and Lajos Soukup\medskip\\}
\date{May 26, 1994}
\newif\iftesting
\newif\ifcommented 
\renewcommand{\baselinestretch}{1.17}
\iftesting
\let\Label\label%
\def\label#1{\marginpar{{\tiny #1}}\Label{#1}}%
\fi
\ifcommented
\fi
{\ifcommented\end{footnotesize}\medskip\\\fi}

\newtheorem{Thm}{{\bf Theorem}}[section]
\newtheorem{Cor}[Thm]{{\bf Corollary}}
\newtheorem{Lemma}[Thm]{{\bf Lemma}}
\newtheorem{Prop}[Thm]{{\bf Proposition}}
\newtheorem{Claim}{{\bf Claim}}[Thm]
\newtheorem{Subclaim}{{\bf Subclaim}}[Claim]

\newcommand{\Thmof}[1]{{Theorem \ref{#1}}}
\newcommand{\Corof}[1]{{Corollary \ref{#1}}}
\newcommand{\Propof}[1]{{Proposition \ref{#1}}}
\newcommand{\Lemmaof}[1]{{Lemma \ref{#1}}}

\newcommand{\Thmabove}{{Theorem \number\theThm}}

\newcommand{\Lemmaabove}{{Lemma \number\theThm}}

\newcommand{\prf}{{\bf Proof\ \ }}

\newcommand{\prfofClaim}{\raisebox{-.4ex}{\Large $\vdash$\ \ }}
\newsavebox{\qedbox}\sbox{\qedbox}{
{\unitlength=0.07mm \begin{picture}(40,60)
\put(0,0){\framebox(30,44)[cc]{}}
\put(30,-7){\rule{7\unitlength}{44\unitlength}}
\put(10,-7){\rule{27\unitlength}{7\unitlength}}
\end{picture}}}
\newcommand{\qed}{\mbox{}\hfill\usebox{\qedbox}}
\newcommand{\smallqed}%
{\mbox{}\smallskip\hfill\raisebox{-.4ex}{\Large $\dashv$}\\}
\newcommand{\qedof}[1]%
{\mbox{} \hspace*{\fill}{\usebox{\qedbox}{~(#1)}}\smallskip\\%
\mbox{}}
\newcommand{\Qedof}[1]%
{\mbox{} \hspace*{\fill}{\usebox{\qedbox}%
{~(#1~\number\theThm)}}\smallskip\\\mbox{}}

\newcommand{\qedofCor}{\Qedof{Corollary}}
\newcommand{\qedofProp}{\Qedof{Proposition}}
\newcommand{\qedofLemma}{\Qedof{Lemma}}

\newcommand{\qedofClaim}%
{\mbox{}\hfill\raisebox{-.4ex}{\Large $\dashv$ }\nolinebreak%
\mbox{~(Claim~\number\theClaim)}\smallskip\\}
\newcommand{\qedofSubclaim}%
{\mbox{}\hfill\raisebox{-.4ex}{\Large $\dashv$ }\nolinebreak%
\mbox{~(Subclaim~\number\theSubclaim)}\smallskip\\}

\newcommand{\assert}[1]{\makebox[6ex]{${\it #1})$}}
\newcommand{\assertof}[1]{${\it #1})$}
\newenvironment{assertion}[1]{\begin{trivlist}
\item[]\hspace*{\parindent}{#1}\hspace{\labelsep}%
\dimen255=\textwidth
\setbox255=\hbox{\hspace*{\parindent}{#1}\hspace{\labelsep}}
\advance\dimen255 by -1\wd255
\advance\dimen255 by -1ex
\begin{minipage}[t]{\dimen255}}%
{\end{minipage}\end{trivlist}}
\newcommand{\implies}{$\,\Rightarrow\,$}

{$\hfill\mbox{}\\[\belowdisplayskip]}
\newcommand{\restr}%
{{\hspace{0.1ex}|\hspace{-0.02ex}{\grave{}}\hspace{0.8ex}}}
\newcommand{\subsetnoneq}%
{\mathrel{\raisebox{-0.8ex}{$\stackrel{\subset}{\scriptstyle\,\not=\,}$}}}

\newcommand{\cardof}[1]{{\mid{#1}\mid}}
\newcommand{\generatedby}[1]{{\langle#1\rangle}}

\newcommand{\setof}[2]{\{\,#1\,:\,#2\,\}}
\newcommand{\smallsetof}[1]{\{\,#1\,\}}

\newcommand{\continuum}{2^{\aleph_0}}

\newcommand{\integers}{Z\hspace*{-1.12ex}Z}
\newcommand{\unitint}{\mbox{$\rm I_{\!\!}I$}}

\newcommand{\forces}[2]{\,\|\hspace{-.35ex}\mbox{\sf--}_{\,#1\,}%
\mbox{\rm``}\,#2\,\mbox{\rm''}}
\let\Models\models
\def\models#1{\Models\mbox{\rm``}\,#1\,\mbox{\rm''}}
\newcommand{\xmbox}[1]{ ${\rm #1}$ }
\newcommand{\st}{such that}
\newcommand{\wolog}{without loss of generality}
\newcommand{\Wolog}{Without loss of generality}
\newcommand{\wrt}{with respect to}

\newcommand{\Ba}{Boolean algebra}
\newcommand{\Bas}{Boolean algebras}
\newcommand{\cBa}{complete Boolean algebra}
\newcommand{\cBas}{complete Boolean algebras}

\newcommand{\po}{partial ordering}
\newcommand{\pos}{partial orderings}

\newcommand{\Fr}{{\rm Fr}\,}

\newcommand{\pcf}{\mathop{\rm pcf}}

\newcommand{\Fn}{{\rm Fn}}
\newcommand{\dom}{{\rm dom}}

\newcommand{\cohenalg}[1]%
{\mbox{$\,\raisebox{0.05ex}{\small$\wr$}\!\!_{_{\!\!}}\mbox{\rm C}_{#1}$}}
\newcommand{\cantorsp}[1]{\mbox{}^{#1^{\mbox{}\!}}2}
\newcommand{\fnsp}[2]{\mbox{}^{#1^{\mbox{}\!}}#2}
\newcommand{\mapping}[3]{#1:#2\rightarrow #3}

\newcommand{\ZFC}{{\rm ZFC}}
\newcommand{\CH}{{\rm CH}}
\newcommand{\MA}{{\rm MA}}

\newcommand{\stick}%
{\mbox{{\hspace{0.4ex}$\mbox{\raisebox{1ex}{$\bullet$}}%
\hspace*{-0.94ex}|$\hspace{0.6ex}}}}

\newcommand{\calC}{{\cal C}}
\newcommand{\calD}{{\cal D}}
\newcommand{\calF}{{\cal F}}

\newcommand{\calP}{{\cal P}}

\newcommand{\dotf}{{\dot{f}}}

\newcommand{\undercirc}[1]{%
\setbox255=\hbox{$#1$}
\dimen255=0.9ex
\advance\dimen255 by \dp255
\rlap{\hbox to\wd255{\hss{\lower\dimen255\hbox{$\scriptstyle\circ$}}\hss}}
\box255
}
\newcommand{\sundercirc}[1]{%
\setbox255=\hbox{$\scriptstyle#1$}
\dimen255=0.8ex
\advance\dimen255 by \dp255
\rlap{\hbox to\wd255{\hss{\lower\dimen255
\hbox{$\scriptscriptstyle\circ$}}\hss}}
\box255
}
\newcommand{\multiundercirc}[1]{%
\mathchoice{\undercirc{#1}}{\undercirc{#1}}{\sundercirc{#1}}{\sundercirc{#1}}}
\newcommand{\acirc}{\multiundercirc{a}}
\begin{document}
\maketitle
\begin{abstract}
We show that a theorem of Leonid B.\ Shapiro which was proved under \MA, is 
actually independent from \ZFC. We also give a direct proof of the Boolean 
algebra version of the theorem under \MA({\it Cohen}). 
\end{abstract}
\renewcommand{\thefootnote}{ }
\footnotetext{The first author would like to thank Professor L.B.\ 
Shapiro for calling his attention to \cite{shapiro}, and 
Professor K.\ Eda for informing him of an example in \cite{eda}. 
The second author is partially supported by the Deutsche 
Forschungs\-gemeinschaft(DFG) grant Ko 490/7--1. 
He also gratefully acknowledges partial support by the Edmund Landau Center 
for research in Mathematical Analysis, supported by the Minerva 
Foundation (Germany). 
The present paper is the second author's Publication No.\ 543. 
The third author is partially supported by the Hungarian National 
Foundation for Scientific Research grant No.\ 1908 and the Deutsche 
Forschungs\-gemeinschaft(DFG) grant Ko 490/7--1.} 
\section{Introduction}
\ifcommented
{\renewcommand{\thefootnote}{\fnsymbol{footnote}}
\footnotetext[2]{the first author visited the Department of Mathematics of 
Silesian University, 
Katowice, Poland at the beginning of December 1993. At the same 
time Professor Leonid Shapiro was also a guest of the university and in 
this way he could learn \Thmof{shapiro} from Professor Shapiro himself. 
He would like to thank Professor Alexander B\l aszczyk 
for arranging this meeting and for his hospitality during his stay in 
Katowice. }
\fi
\ifcommented
\footnotetext[3]{This note grew out of a seminar talk given by the first 
author in Berlin in January 1994. 
}} \fi
\noindent
L.B.\ Shapiro \cite{shapiro} recently proved the following theorem: 
\begin{Thm}{\rm (L.B.\ Shapiro) (\MA({\small\it Cohen}))}\label{shapiro}
For any compact Hausdorff space $X$ of weight $<\continuum$ and 
$\aleph_0\leq\tau<\continuum$ the following assertions are equivalent:
\medskip\\
\assert{i} There exists a continuous surjection from $X$ onto 
$\fnsp{\tau}{\unitint}$;\\ 
\assert{ii} There exists a continuous injection from 
$\cantorsp{\tau}$ into $X$;\\ 
\assert{iii} There exists a closed subset $Y\subseteq X$ \st\ 
$\chi(y,Y)\geq\tau$ for every $y\in Y$.
\end{Thm}%
The original proof of \Thmabove\ by L.B.\ Shapiro in \cite{shapiro} was 
formulated under \MA. However practically the same proof still works when 
merely \MA({\it\small Cohen}) is assumed
where \MA({\it\small Cohen}) stands for Martin's Axiom restricted to the \pos\ of 
the form $\Fn(\kappa,2)$. 

A part of the theorem above can be translated into the 
language of \Bas:
\begin{Cor}\label{shapiro2}{\rm(Boolean algebra version of Shapiro's 
theorem)} 
{\rm(\MA({\small\it Cohen}))} 
For any infinite \Ba\ $B$ of cardinality $<\continuum$ and any infinite 
$\tau$, the following are 
equivalent:\medskip\\
\assert{i'} There exists an injective Boolean mapping from 
$\Fr\tau$ into $B$;\\ 
\assert{ii'} There exists a surjective Boolean mapping from $B$ onto $\Fr\tau$.
\end{Cor}

The implication from \assertof{ii'} to \assertof{i'} as well as the 
implication from \assertof{ii} to \assertof{i} can be proved already in \ZFC. 
For the proof of \assertof{ii} from \assertof{i}, let 
$\mapping{g}{\cantorsp{\tau}}{X}$ be a continuous injection. 
Note that $g[\cantorsp{\tau}]$ is a closed subset of $X$.
For any fixed $y_0\in\cantorsp{\tau}$ let $\mapping{f'}{X}{\cantorsp{\tau}}$ 
be defined by
\[ f'(x)=\left\{
\begin{array}{ll}
g^{-1}(x)&\mbox{\ \ ; if }x\in g[\cantorsp{\tau}],\\
y_0&\mbox{\ \ ; otherwise.}
\end{array}
\right.
\]\noindent
Then $f'$ is a continuous surjection from $X$ onto $\cantorsp{\tau}$. 
Let $f''$ be a continuous surjection from 
$\cantorsp{\tau}$ to $\fnsp{\tau}{\unitint}$.\label{jousha}
E.g.\ let 
$h:\cantorsp{\omega}\rightarrow \unitint$ be the continuous surjection 
defined by $u\mapsto$ the real represented by the 
binary expression $0.u(0)u(1)u(2)\cdots$. 
$\mapping{h^\kappa}%
{\fnsp{\kappa}{(\cantorsp{\omega})}}{\fnsp{\kappa}{\unitint}}$
is then a continuous surjection. Since $\fnsp{\kappa}{(\cantorsp{\omega})}$ 
is homeomorphic to $\cantorsp{\kappa}$ we can find a continuous surjection 
$f''$ from $\cantorsp{\tau}$ onto $\fnsp{\tau}{\unitint}$ corresponding to 
$h^\kappa$. 
The mapping $g=f''\circ f'$ is then as desired.
In the next section we shall give a direct proof 
of \assertof{i'}\implies\assertof{ii'}. 
For \assertof{iii}\implies\assertof{i} we need some deep results by 
Shapiro on dyadic compactum (see \cite{shapiro}).

The equivalence of the assertions \assertof{i'} and \assertof{ii'} 
above is not true in general for \Bas\ of cardinality 
$\geq \continuum$: For any $\sigma$-complete \Ba\ $B$ and any infinite 
$\kappa$, there 
exits no surjective Boolean mapping $\mapping{f}{B}{\Fr\kappa}$ 
(see \Lemmaof{no-surj} below). 
Hence e.g.\ for \Ba\ $B=\overline{\Fr\omega}$ we have that 
$\cardof{B}=\continuum$; $\Fr \continuum$ is embeddable into $B$ (by 
Balcar-Fra\v nek-Theorem, see \cite{balcar-franek}) but 
there exists no surjective Boolean mapping from $B$ onto $\Fr \continuum$.
The non-existence of surjective Boolean mapping from a 
$\sigma$-complete \Ba\ in the ground model onto $\Fr\tau$ is preserved 
in a generic extension by a \po\ of cardinality $<\tau$ though $B$ may be 
no more $\sigma$-complete in such a generic extension:
\begin{Lemma}\label{no-surj}
Let $B$ be a $\sigma$-complete \Ba\ and $P$ a \po. For any 
$\kappa>\cardof{P}$ we have that
\[\forces{P}{\mbox{there exists no surjective Boolean mapping from }
B\mbox{ onto }\Fr\kappa}. 
\]\noindent
\end{Lemma}
\prf
Suppose that there would be a $P$-name $\dotf$ \st\ 
\[ \forces{P}{\mapping{\dotf}{B}{\Fr\kappa}
\mbox{ is a surjective Boolean mapping}}.
\]\noindent
For each $p\in P$ let 
\[ B_p=\setof{b\in B}{p\forces{P}{\dotf(b)=c\mbox{ for some }c\in\Fr\kappa}}
\]\noindent
and
\[ C_p=\setof{c\in\Fr\kappa}{p\forces{P}{\dotf(b)=c\mbox{ for some }b\in B}}.
\]\noindent
Then $B_p$ and $C_p$ are subalgebras of $B$ and $\Fr\kappa$ respectively. 
Since $\bigcup_{p\in P}C_p=\Fr\kappa$ and $\kappa>\cardof{P}$ there exists 
some $p\in P$ \st\ $C_p$ is infinite. Let $c_n$, $n<\omega$ be 
pairwise disjoint positive elements of $C_p$. By the definition of 
$B_p$ and $C_p$, there exits pairwise disjoint positive elements $b_n$, 
$n<\omega$ of $B_p$ \st\ $p\forces{P}{\dotf(b_n)=c_n}$ holds for every 
$n<\omega$. 
Let $X\subseteq\omega$ be \st\ there exists no $c\in\Fr\kappa$ \st\ 
$c\cdot c_n=c_n$ holds for all $n\in X$ and $c\cdot c_n=0$ for all 
$n<\omega\setminus X$. 
Let $d=\Sigma^B_{n\in X}b_n$.
Then for any $q\leq p$ there can be no $c\in\Fr\kappa$ \st\ 
$q\forces{P}{\dotf(d)=c}$. This is a contradiction. 
\qedofLemma
The lemma above together with \Corof{shapiro2} yields the following:
\begin{Prop}
Let $B$ be a \cBa\ with $\cardof{B}=\tau\geq\aleph_0$. Then 
\[ \forces{\Fn(\kappa,2)}%
{\mbox{there exists no surjective Boolean mapping from }
B\mbox{ onto }\Fr\tau}
\]\noindent
holds if and only if $\kappa<\tau$.
\end{Prop}
\prf
If $\kappa<\tau$ then $\cardof{\Fn(\kappa,2)}=\kappa<\tau$. Hence 
by \Lemmaof{no-surj}, 
\[\forces{\Fn(\kappa,2)}%
{\mbox{there exists no surjective Boolean mapping from }
B\mbox{ onto }\Fr\tau}
\]\noindent
holds. 

Suppose now that $\kappa\geq\tau$. 
Then as in the proof of \Propof{shapiro3}, we can show that
\[ \forces{\Fn(\kappa,2)}%
{\mbox{there exists a surjective Boolean mapping from }
B\mbox{ onto }\Fr\tau}
\]\noindent
holds.
\qedofProp

Now, (\stick)\ (read ``stick'', see \cite{broverman-ginsburg-kunen-tall}) 
is the following principle:
\begin{assertion}{(\stick):}
 There exists a sequence $(x_\alpha)_{\alpha<\omega_1}$ of 
countable subsets of $\omega_1$ \st\ for any $y\in[\omega_1]^{\aleph_1}$ 
there exists $\alpha<\omega_1$ \st\ $x_\alpha\subseteq y$.
\end{assertion}
Clearly (\stick) follows from \CH. Another combinatorial principle 
$(\clubsuit)$, a strengthning of (\stick), is 
introduced in Ostaszewski \cite{ostaszewski}. Let 
$Lim(\omega_1)=\setof{\gamma<\omega_1}{\gamma\mbox{ is a limit}}$.
\begin{assertion}{$(\clubsuit)$:} There exists a sequence 
$(x_\gamma)_{\gamma\in Lim(\omega_1)}$ of countable 
subsets of $\omega_1$ \st\ for every $\gamma\in Lim(\omega_1)$, $x_\gamma$ is 
a cofinal subset of $\gamma$, $otp(x_\gamma)=\omega$ and for every 
$X\in[\omega_1]^{\aleph_1}$ there is 
$\gamma\in Lim(\omega_1)$ \st\ $x_\gamma\subseteq X$.
\end{assertion}
Clearly (\stick) follows from $(\clubsuit)$. 
Unlike (\stick), $(\clubsuit)$ does not follow from \CH, since  
$(\clubsuit)$ $+$ \CH\ is equivalent with $\Diamond$ (K.\ Devlin, 
see \cite{ostaszewski}). For more about the combinatorial principles (\stick) 
and $(\clubsuit)$, and independence results connected with them, 
see \cite{FShS544}. 

\MA({\small\it countable}) --- Martin's axiom restricted to 
countable \pos\ --- and \MA({\small\it Cohen}) both add a lot of Cohen 
reals over 
any small model of (a sufficiently large finite subset of) \ZFC\ and in
many cases where this property is needed, \MA({\small\it countable}) is just 
enough. Hence it seems to be quite natural to ask if these axioms are 
perhaps equivalent. 
However they are not. 
I.\ Juh\'asz proved in an  unpublished note that 
$\neg\CH$ $+$ \MA({\small\it countable}) 
$+$ $(\clubsuit)$ is consistent (two other constructions of models of 
$\neg\CH$ $+$ \MA({\small\it countable}) 
$+$ $(\clubsuit)$ are to be found in \cite{komjath} and \cite{FShS544}.). 
On the other hand, it is easy 
to see that the negation of $\MA(\Fn(\aleph_1,2))$ follows from 
$\neg\CH$ $+$ $(\clubsuit)$: 
using $(\stick)$ we can obtain a \Ba\ $B$ of cardinality 
$\aleph_1$ \st\ $\Fr\omega_1$ is embeddable into $B$ but there is no 
surjection from $B$ onto $\Fr\omega_1$ (see \Thmof{non-shapiro}). 
By \Propof{shapiro3}, this shows that 
$m_{\Fn(\aleph_1,2)}=\aleph_1<\continuum$. 
It follows also that the assertions of \Thmof{shapiro} 
and \Corof{shapiro2} are independent from \ZFC\ and \MA({\small\it countable}) 
is not 
enough to prove them. 

\Corof{shapiro2} for other variety than \Bas\ can be simply false. 
E.g., this is the case in the variety of abelian groups:
in \cite{eda}, an $\aleph_1$-free abelian group $G$ in $\aleph_1$ is 
constructed (in \ZFC) which contains uncountable free subgroup but 
$Hom(G,Z)=0$. 

\section{A proof of the Boolean algebra version of the theorem}
In this section we shall prove \Corof{shapiro2}. More precisely we 
prove the following \Propof{shapiro3}. For any class $\calC$ of \pos\ Let
\[ 
\begin{array}{@{}l@{}l}
m_{\calC}=\min\setof{\cardof{\calD}}{&\calD\mbox{ is a family of dense 
subsets of }P\mbox{ for some }P\in\calC\\
&\mbox{\st\ there exists no }\calD\mbox{-generic filter over }P
}
\end{array}
\]\noindent
If $\calC$ is a singleton $\smallsetof{P}$, we shall write simply $m_P$ in 
place of $m_{\smallsetof{P}}$. 
Let us say that two \pos\ $P$, $Q$ are coabsolute when their 
completions are isomorphic. It is easy to see that for any class 
$\calC$ of \pos\ $m_\calC=m_{\tilde\calC}$ where 
$\tilde\calC=\setof{Q}{Q\xmbox{ is coabsolute with some }P\in\calC}$. 
If the class $\calC$ is introduced by a property $\calP$ of \Bas, we also 
write $m_\calP$ in place of $m_\calD$. 
We also write 
$m_{\mbox{\footnotesize\it countable}}=
m_{\setof{P}{P\mbox{\footnotesize\ is countable}}}$ and 
$m_{\mbox{\footnotesize\it Cohen}}=
m_{\setof{P}{P=\Fn(\kappa,2)\mbox{\footnotesize\ for some }\kappa}}$. 
Hence \MA({\small\it Cohen}) (\MA({\small\it countable}), \MA\ etc.\ 
respectively) 
holds if and only if 
$m_{\mbox{\footnotesize\it Cohen}}=\continuum$ 
($m_{\mbox{\footnotesize\it countable}}=\continuum$, $m_{ccc}=\continuum$ 
etc.\ respectively)   
and we have $m_{ccc}\leq m_{\mbox{\footnotesize\it Cohen}}\leq 
m_{\mbox{\footnotesize\it countable}}$. 
\begin{Prop}\label{shapiro3}
Let $B$ be a \Ba\ containing $\Fr\kappa$ as a subalgebra. 
If $\cardof{B}<m_{\Fn(\kappa,2)}$, 
then there exists a surjective Boolean mapping from $B$ onto $\Fr\kappa$.
\end{Prop}
\prf
By Sikorski's theorem, there is a Boolean mapping from $B$ to 
$\overline{\Fr\kappa}$ --- the completion of $\Fr\kappa$, extending the 
inverse of the canonical embedding of $\Fr\kappa$ into $B$. 
Hence \wolog\ we may assume 
that $B$ is a subalgebra of $\overline{\Fr\kappa}$. 
Now let $P=\Fn(\kappa,3)$. Note that $P$ is coabsolute with $\Fn(\kappa,2)$. 
We shall define a family $\calD$ of dense subsets of 
$P$ \st\ $\cardof{D}< m_{\Fn(\kappa,2)}$ so that among other things (see 
below), for $\calD$-generic set $G$, $g=\bigcup G$ will be a function from 
$\kappa$ to $3$ and $X=\setof{\alpha<\kappa}{g(\alpha)=2}$ will be 
of cardinality $\kappa$. 
Then we let $f$ be the function on $\kappa$ defined by: 
\[ f(\alpha)=\left\{
\begin{array}{ll}
	0_B\quad&\mbox{; if }g(\alpha)=0,\\
	1_B\quad&\mbox{; if }g(\alpha)=1,\\
	\alpha&\mbox{; otherwise.}
\end{array}
\right.
\]\noindent
Let $\bar f$ be the Boolean mapping from $\Fr\alpha$ to $\Fr X$ 
generated by $f$. 

Now we are done, if we can show that $\bar f$ extends to a Boolean mapping 
$\tilde f$ from $B$ onto $\Fr X$. But by the following \Lemmaof{LemmaX}, 
we can choose $\calD$ appropriate for this purpose.

For $p\in P$, let $B_p=\Fr\dom(p)$ (hence $B_p\leq B$) and 
$\mapping{f_p}{B_p}{\Fr(p^{-1}[\{2\}])}$ be the Boolean mapping generated 
by the mapping $f^0_p$ on $\dom(p)$ defined by:
\[ f^0_p(\alpha)=\left\{
\begin{array}{ll}
	0_{B}\quad&\mbox{; if }p(\alpha)=0,\\
	1_{B}\quad&\mbox{; if }p(\alpha)=1,\\
	\alpha&\mbox{; otherwise}.
\end{array}
\right.
\]\noindent
\begin{Lemma}\label{LemmaX}
For any $b\in B$ and $p\in P$ there exists $q\leq p$ and $b_1$, 
$b_2\in B_q$ \st\ $b_1\leq b$, $b_2\leq -b$ and $f_q(b_1)+f_q(b_2)=1$ {\rm(}
i.e, $q$ ``forces'' $\tilde f(b)=f_q(b_1)$~{\rm)}.
\end{Lemma}
For the proof of the \Lemmaabove\ we use the following Lemma whose proof is 
left to the reader:
\begin{Lemma}
Let $b\in\overline{\Fr\kappa}$ and let 
$Y\subseteq\kappa$ be a countable set \st\ $b\in\overline{\Fr Y}$ holds. 
Let 
$Y=\setof{\alpha_n}{n<\omega}$. Then there exist an increasing 
sequence $(l_n)_{n<\omega}$ with $l_n<\omega$ for $n<\omega$ 
and a sequcence $(i_n)_{n<\omega}$ with 
$i_n\in\fnsp{l_n}{\smallsetof{-1,1}}$ for 
$n<\omega$ \st, letting $p_n=\Sigma_{k<l_n}i_n(k)\cdot\alpha_k$ for 
$n<\omega$, \smallskip\\
\assert{i} either $p_n\leq b$ or $p_n\leq -b$ and \\
\assert{ii} $\Sigma_{n<\omega}p_n=1$.
\smallskip\\
In particular we have $b=\Sigma\setof{p_n}{n<\omega,\, p_n\leq b}$.\qed
\end{Lemma}
{\bf Proof of \Lemmaof{LemmaX}}\quad 
Let $Y=\setof{\alpha_n}{n<\omega}$, $(l_n)_{n<\omega}$, 
$(i_n)_{n<\omega}$ and $p_n$, $n<\omega$ be as in \Lemmaabove\ for our 
$b\in B$. 
\Wolog\ we may assume that $\dom(p)\cap Y=\setof{\alpha_n}{n<k}$ for 
some $k<\omega$. Let $\fnsp{k}{\smallsetof{-1,1}}=\setof{\tau_m}{m<2^k}$. 
By induction we can take $n_m<\omega$ for $m<2^k$ \st\medskip\\
\assert{a} $i_{n_m}$ is compatible (as an element of 
$\Fn(Y,\smallsetof{-1,1})$) 
with $\tau_m$ and\smallskip\\
\assert{b}
$\setof{i_{n_m}\restr(\dom(i_{n_m})\setminus k)}{m<2^k}$ is pairwise 
compatible.\medskip\\ 
Let $\tilde n=\max\setof{n_m}{m<2^k}$, $\tilde l=l_{\tilde n}$ and 
$\tilde i=\bigcup\setof{i_{n_m}\restr(\dom(i_{n_m})\setminus k)}{m<2^k}$.
Let $q\leq p$ be \st\ $\dom(q)=\dom(p)%
\cup\smallsetof{\alpha_k,\ldots,\alpha_{\tilde l-1}}$,
$q\restr \dom(p)=p$ and
\[ q(\alpha_m)=\left\{
\begin{array}{ll}
1\quad&\mbox{; if }\tilde i(\alpha_m)=1,\\
0\quad&\mbox{; if }\tilde i(\alpha_m)=-1.
\end{array}
\right.
\]\noindent
Then $q$ as above together with 
$b_1=\Sigma\setof{p_n}{n<\tilde n,\,p_n\leq b}$	and
$b_2=\Sigma\setof{p_n}{n<\tilde n,\,p_n\leq -b}$ 
is as desired.\qedof{\Lemmaof{LemmaX}}
Now by the lemma above 
\[ 
\begin{array}{@{}l@{}l}
\calD\,=\,
&\setof{\setof{p\in P}{\alpha\in\dom(p)}}{\alpha<\kappa}\\
&\cup\setof{\setof{p\in P}{\exists\beta>\alpha\,\, p(\beta)=2}}%
{\alpha<\kappa}\\
&\cup\setof{\setof{q\in P}{f_q(b_1)+f_q(b_2)=1
\mbox{ for some }b_1\leq b,\,b_2\leq -b}}{b\in B}
\end{array}
\]\noindent
is a family of dense subsets of $P$. Clearly the mapping $\bar{f}$ defined 
as above \wrt\ this $\calD$ can be extended to a Boolean mapping $\tilde f$ 
from $B$ onto $\Fr X$. 
\qedof{\Propof{shapiro3}}
\section{Pcf and the theorem of Shapiro}\label{pcf}
\begin{Prop}\label{oplus}
Assume that
\begin{assertion}{$\oplus_{\mu,\kappa,\lambda}$}
for any $\calF\subseteq[\lambda]^{\aleph_0}$ with $\cardof{\calF}<\mu$, 
there is $Y\in[\lambda]^{\kappa}$ \st\ $a\cap Y$ is finite for all 
$a\in\calF$. 
\end{assertion}
Then, for any \Ba\ $B$ of cardinality $<\mu$, if $\Fr\lambda$ is embeddable 
into $B$ then there is a surjective Boolean mapping from $B$ onto 
$\Fr\kappa$. 
\end{Prop}
\prf
As in the proof of \Propof{shapiro3}, we may assume \wolog\ that 
$\Fr\lambda\leq B\leq\overline{\Fr\lambda}$ holds. Let 
$\cardof{B}=i^*$ ($<\mu$) and let $(y_i)_{i<i^*}$ be an enumeration of $B$. 
Let 
$y_i=\sum_{n<\omega}\tau^n_i(\alpha(i,n,0),\ldots,\alpha(i,n,m_{i,n}))$ 
where $\tau^n_i$ is a Boolean term with $m_{i,m}+1$ variables and 
$\alpha(i,n,0),\ldots,\alpha(i,n,m_{i,n})<\lambda$ for $i<i^*$ and 
$n<\omega$. For $i<i^*$, let 
$w_i=\setof{\alpha(i,n,l)}{n<\omega,\,l\leq m_{i,n}}$. By the assumption, 
there exists a $Y\in[\lambda]^\kappa$ \st\ $w_i\cap Y$ is finite for everly 
$i<i^*$. Let $\mapping{g}{B}{\Fr Y}$ be defined by 
\[ g(y_i)=\sum_{n<\omega}\tau^n_i(\alpha^*(i,n,0),
\ldots,\alpha^*(i,n,m_{i,n}))
\]\noindent
where
\[ \alpha^*(i,n,l)=\left\{%
\begin{array}{ll}
\alpha(i,n,l)\qquad&;\mbox{ if }\alpha(i,n,l)\in Y\\
0_B&;\mbox{ otherwise}.
\end{array}\right.
\]\noindent
The function $g$ is well-defined since, for each $i<\omega$, 
$\tau^n_i(\alpha^*(i,n,0), 
\ldots,\alpha^*(i,n,m_{i,n}))$ is an element of $\Fr(w_i\cap Y)$ and 
$\Fr(w_i\cap Y)$ is 
finite. Clearly this $g$ is as desired. 
\qedofProp
\begin{Cor}
For any \Ba\ of cardinality $<{\bf a}$ (where $\bf a$ is the minimal 
cardinality of a maximal almost disjoint family in $[\omega]^{\aleph_0}$), 
if\/ 
$\Fr\omega$ is embeddable into $B$ then there is a surjection from $B$ onto 
$\Fr\omega$. 
\end{Cor}
\prf By \Propof{oplus} for $\oplus_{{\bf a},\aleph_0,\aleph_0}$.\qedofCor

\begin{Thm}\label{pcf-shapiro}
Assume that
\begin{assertion}{$(*)_{\mu,\lambda,\kappa}$}
there are $\acirc_i\in[Reg\cap(\lambda^+\setminus\kappa^+)]^{<\aleph_0}$ for 
$i<\kappa$ \st\ for every $a\in[\kappa]^{\aleph_0}$, 
$\max\,\pcf(\bigcup_{i\in a}\acirc_i)\geq\mu$ holds.
\end{assertion}
Then for any \Ba\ $B$ of cardinality $<\mu$, if $\Fr\kappa$ is embeddable 
into $B$ then there is a surjective Boolean mapping $g$ from $B$ onto 
$\Fr\kappa$.  
\end{Thm}
(\,For more about $(*)_{\mu,\lambda,\kappa}$ see {\rm\cite{sh513}}. For pcf 
theory in general, the reader may consult \cite{sh-g}.) 
The theorem follows from \Propof{oplus} and the following:
\begin{Lemma}
Assume that $(*)_{\mu,\lambda,\kappa}$ (as in \Thmof{pcf-shapiro}) holds. 
Then $\oplus_{\mu,\kappa,\kappa}$ holds.
\end{Lemma}
\prf
Since $\max\,\pcf$ is always regular, we may assume that $\mu$ is regular. 
Let $\acirc=\bigcup_{i<\kappa}\acirc_i$. 
In place of $[\kappa]^{\aleph_0}$, we consider $[Z]^{\aleph_0}$ for 
$Z=\mathop{\stackrel{\cdot}{\bigcup}}_{i<\kappa}Z_i$ where 
$Z_i=\smallsetof{i}\times\prod\acirc_i$. Hence we assume that 
$\calF\subseteq[Z]^{\aleph_0}$ and $\cardof{\calF}<\mu$. 

For each $a\in\calF$, let $g_a\in\prod\acirc$ be defined by 
\[ g_a(\theta)=\sup\setof{\eta(\theta)}{\eta\in a,\,\theta\in\dom(\eta)}
\]\noindent
for each $\theta\in\acirc$, where we put $\sup\emptyset =0$. 
Since $\prod\acirc/J_{<\mu}[\acirc]$ is $\mu$-directed and 
$\cardof{\calF}<\mu$, there is 
$f^*\in\prod\acirc$ \st\ $g_a<_{J_{<\mu}[\acirc]}f^*$ holds for all 
$a\in\calF$. 
For $i<\kappa$, let $z_i=\smallsetof{(0,i)}%
\cup(f^*\restr\acirc_i)$. Then $z_i\in Z_i$ for $i<\kappa$. We 
show that $Y=\setof{z_i}{i<\kappa}$ is as required. Suppose not. Then 
$Y\cap a$ would be infinite for some $a\in\calF$. By the assumption, it 
follows that 
$\bigcup_{z_i\in Y\cap a}\acirc_i\not\in J_{<\mu}[\acirc]$. But for 
$z_i\in Y\cap a$ we have 
$\smallsetof{(0,i)}\cup(f^*\restr\acirc_i)\in a$. It follows that for 
$\theta\in\acirc_i$ we have $f^*(\theta)\leq g_a(\theta)$. This is 
a contradiction to $g_a<_{J_{<\mu}[\acirc]}f^*$.
\qedofLemma
\section{Independence of the theorem of Shapiro}\label{independence}
The principle (\stick) suggests the following cardinal invariant \stick:
\[ \stick= \min
\setof{\cardof{X}}{X\subseteq[\omega_1]^{\aleph_0},\,
\forall y\in[\omega_1]^{\aleph_1}\:\exists x\in X\;x\subseteq y}.
\]\noindent
Clearly $\aleph_1\leq \stick\leq\continuum$ and (\stick) holds if and only if 
$\stick=\aleph_1$. We can also consider the following variants of $\stick$:
\[ 
\begin{array}{@{}l@{}l}
\stick'=\min
\setof{\kappa}{&\kappa\geq\aleph_1,\,\mbox{there is an }
X\subseteq[\kappa]^{\aleph_0}\\
&\mbox{\st\ }\cardof{X}=\kappa\mbox{ and }
\forall y\in[\kappa]^{\aleph_1}\:\exists x\in X\;x\subseteq y},
\end{array}
\]\noindent
\[ 
\begin{array}{@{}l@{}l}
\stick''=\min
\setof{\kappa}{&\kappa\geq\aleph_1,\,\mbox{there is an }
X\subseteq[\kappa]^{\aleph_0}\\
&\mbox{\st\ }\cardof{X}=\kappa\mbox{ and }
\forall y\in[\kappa]^{\kappa}\:\exists x\in X\;x\subseteq y}.
\end{array}
\]\noindent
We have $\aleph_1\leq \stick''\leq \stick'\leq\continuum$ and 
(\stick) holds if and only if $\stick=\stick'=\stick''=\aleph_1$ holds.

It can be easily shown that $\stick\leq\stick'$ holds. Moreover if 
$\stick<\aleph_{\omega_1}$, then $\stick=\stick'$ holds. The question, if 
$\stick<\stick'$ is consistent, is connected with some very fundamental 
unsolved problems on cardinal arithmetics while we can show that 
$\stick''<\stick$ is consistent. For more, see \cite{FShS544} 
and \cite{sh513}.

\begin{Prop}\label{Ba-in-st'}
There exists a \Ba\ $B$ \st\ $\cardof{B}=\stick'$, 
$\Fr \stick'$ is embeddable into 
$B$ but there is no surjective Boolean mapping from $B$ onto $\Fr\omega_1$. 
\end{Prop}
\prf
Let $\mapping{\Phi}{\kappa}{\kappa}$; $\alpha\mapsto\xi_\alpha$ be the 
continuously increasing function defined inductively by 
$\xi_0=\omega$ and $\xi_{\alpha+1}=\xi_\alpha+\cardof{\xi_\alpha}$. 
Let $\kappa=\stick'$ and let $X\subseteq[\kappa\times\Fr\omega_1]^{\aleph_0}$ 
be \st\ $\cardof{X}=\kappa$, $\omega\times\Fr\omega\in X$ and 
$\forall y\in[\kappa\times\Fr\omega_1]^{\aleph_1}\:
\exists x\in X\;x\subseteq y$ holds. 
Let $(x_\alpha)_{\alpha<\kappa}$ be an enumeration of 
$X$ \st\ $x_\alpha\subseteq\xi_\alpha\times\Fr\omega_1$ for all 
$\alpha<\kappa$. 

Now let $(B_\alpha)_{\alpha<\kappa}$ be a continuously increasing sequence 
of \Bas\ \st\ for all $\alpha<\kappa$
\begin{itemize}
\setlength{\itemsep}{\medskipamount}
\item[{\it 1})] the underlying set of $B_\alpha$ is $\xi_\alpha$;
\item[{\it 2})] there exits a $b_\alpha\in B_{\alpha+1}$ \st\ $b_\alpha$ is 
free 
over $B_\alpha$;
\item[{\it 3})] if $x_\alpha$ generates a Boolean mapping $f_\alpha$ from a 
subalgebra of $B_\alpha$ onto an infinite subalgebra of $\Fr\omega_1$ then 
$B_{\alpha+1}$ contains an element $c_\alpha$ of the form 
$\Sigma_{n\in Z_\alpha}^{B_{\alpha+1}}b^\alpha_n$ where 
$Z_\alpha\subseteq\omega$, $b^\alpha_n$, $n<\omega$ are pairwise disjoint 
elements in $\dom(f_\alpha)$, $f_\alpha(b^\alpha_n)\not=0$ for all 
$n<\omega$ and there is no 
$d\in\Fr\omega_1$ \st\ $d\cdot f_\alpha(b^\alpha_n)=f(b^\alpha_n)$ for all 
$n\in Z_\alpha$ and $d\cdot f_\alpha(b^\alpha_n)=0$ for all 
$n<\omega\setminus Z_\alpha$ holds.
\end{itemize}
\par
Let $B=\bigcup_{\alpha<\kappa}B_\alpha$. We show that this $B$ is as 
desired. By \assertof{1} the underlying set of $B$ is $\kappa$. 
By \assertof{2} $\setof{b_\alpha}{\alpha<\kappa}$ is an independent subset 
of $B$. Hence $\Fr\kappa$ is embeddable into $B$. 

Suppose now that there would be a surjective Boolean mapping $f$ from $B$ 
onto $\Fr\omega_1$. Then there is a bijection $g\subseteq f$ from a subset 
of $B$ onto $\Fr\omega_1$. Since $g$ is uncountable there is an 
$\alpha<\kappa$ \st\ $x_\alpha\subseteq g$. Since $x_\alpha\subseteq f$, 
$x_\alpha$ satisfies the condition in \assertof{3}. Hence there is 
a $c_\alpha\in B_{\alpha+1}$ \st\ 
$c_\alpha=\Sigma_{n\in Z_\alpha}^{B_{\alpha+1}}b^\alpha_n$ for 
$Z_\alpha$ and $b^\alpha_n$, $n<\omega$ as un \assertof{3}. But then  
$f(c_\alpha)\cdot f_\alpha(b^\alpha_n)=f(b^\alpha_n)$ for all 
$n\in Z_\alpha$ and $f(c_\alpha)\cdot f_\alpha(b^\alpha_n)=0$ for all 
$n<\omega\setminus Z_\alpha$ holds. This is a contradiction to the choice 
of $Z_\alpha$. 
\qedofProp
\begin{Cor}\label{m-Fn-leq-st'}
$m_{\Fn(\omega_1,2)}\leq \stick'$.
\end{Cor}
\prf
By \Propof{shapiro3} and \Propof{Ba-in-st'}.\qedofCor
With almost the same proof as in \Propof{Ba-in-st'} we can also prove the 
following: 
\begin{Prop}\label{Ba-in-st''}
There exists a \Ba\ $B$ \st\ $\cardof{B}=\stick''$, 
$\Fr\stick''$ is embeddable into 
$B$ but there is no surjective Boolean mapping from $B$ onto $\Fr\stick''$. 
\qed
\end{Prop}
Since we have $\stick'=\aleph_1$ under $(\stick)$, 
we obtain the following theorem:
\begin{Thm}\label{non-shapiro}
If $(\stick)$ holds then there exists a \Ba\ $B$ of cardinality 
$\aleph_1$ \st\ $\Fr\omega_1$ is embeddable into $B$ but there is no 
surjection from $B$ onto $\Fr\omega_1$.\\\mbox{}\qed
\end{Thm}
Hence if $\neg$\CH\ and $(\stick)$ holds, by \Thmabove, there exists 
a counter-example to the theorem of Shapiro. 
This shows that we cannot just drop \MA({\small\it Cohen}) 
from \Thmof{shapiro}. 
Since $\MA(\mbox{\small\it countable})$ $+$ $\neg$\CH\ $+$ $(\stick)$ is 
consistent (see e.g.\ \cite{komjath} or \cite{FShS544}), we see that 
\MA({\small\it countable}) is not enough for \Thmof{shapiro}.
\begin{Cor}\label{m-cohen-leq-st''}
$m_{\mbox{\footnotesize\it Cohen}}\leq \stick''$.
\end{Cor}
\prf
By \Propof{shapiro3} and \Propof{Ba-in-st''}.
\qedofCor

If a \Ba\ $B$ is atomless then $\Fr\omega$ can be embdded into $B$. 
By \Propof{shapiro3}, if \MA({\small\it countable}) holds and $B$ is of 
cardinality $<\continuum$, there exists a surjection from $B$ onto 
$\Fr\omega$. Here again we cannot simply drop the assumption of 
\MA({\small\it countable}): 
\begin{Prop}
It is consistent that there is an atomless \Ba\ $B$ of cardinality 
$\aleph_1<\continuum$ \st\ there is no surjective Boolean mapping from $B$ 
onto $\Fr\omega$. 
\end{Prop}
\prf
By \cite[Theorem 5.12]{sh130}, there is a model of \ZFC $+$ $\neg$\CH\ in 
which there is an endo-rigid atomless \Ba\ $B$ of cardinality $\aleph_1$. 
In particular there is no surjection from $B$ onto $\Fr\omega$. 
\qedofProp

Note that, since $(\stick)$ is consistent with 
$\neg$\CH\ and \MA({\small\it countable}),  
$(\stick)$ cannot supply such a \Ba\ as in the proposition above.

\renewcommand{\baselinestretch}{0.9}
\small
\mbox{}
\vfill\mbox{}\\
\mbox{}\qquad\qquad\qquad{\bf Authors' addresses:}\bigskip\\
\mbox{}\qquad\qquad\qquad\qquad
\parbox[t]{12cm}{\it
Institut of Mathematics, Hebrew University of Jerusalem\\
91904 Jerusalem, Israel\smallskip\\
and\smallskip\\
Institut f\"ur Mathematik II, Freie Universit\"at Berlin\\
14195 Berlin, Germany\medskip\\
{\tt fuchino@math.fu-berlin.de}
\bigskip\bigskip\\
Institut of Mathematics, Hebrew University of Jerusalem\\
91904 Jerusalem, Israel\smallskip\\
and\smallskip\\
Department of Mathematics, Rutgers University\\
New Brunswick, NJ 08854, USA\medskip\\
{\tt shelah@math.huji.ac.il}
\bigskip\bigskip\\
Mathematical Institute of the Hungarian Academy of Sciences\medskip\\
{\tt soukup@math-inst.hu}
}

\end{document}